\documentclass[11pt]{amsart}

\usepackage{amsmath}
\usepackage{amssymb}
\usepackage{latexsym}
\usepackage{amscd}
\usepackage[frame,cmtip,arrow,matrix,line,graph,curve]{xy}

\newcommand{\C}{\mathbb{C}}
\newcommand{\PP}{\mathbb{P}}
\newcommand{\m}{\mathfrak{m}}
\newcommand{\Homom}{\mathrm{Hom}}
\newcommand{\Endom}{\mathrm{End}}
\newcommand{\Ker}{\mathrm{Ker}}
\newcommand{\Coker}{\mathrm{Coker}}
\newcommand{\Iden}{\mathrm{Id}}
\newcommand{\rank}{\mathrm{rk}}
\newcommand{\Sym}{\mathrm{Sym}}
\newcommand{\uKer}{\mathrm{\underline{Ker}}}
\newcommand{\uCoker}{\mathrm{\underline{Coker}}}
\newcommand{\Ox}{O_{X}}
\newcommand{\Kx}{K_{X}}
\newcommand{\cI}{{\mathcal I}}
\newcommand{\Spec}{\mathrm{Spec}}
\newcommand{\modulo}{\: \mathrm{mod} \:}
\newcommand{\codim}{\mathrm{codim}}
\newcommand{\Oxx}{{\mathcal O}_{X,x}}
\newcommand{\Span}{\mathrm{Span}}
\newcommand{\ogamma}{\overline{\gamma}}
\newcommand{\Osc}{\mathrm{Osc}}
\newcommand{\U}{{\mathcal U}}
\newcommand{\Cone}{\mathrm{Cone}}

\newtheorem{thm}{Theorem}[section]
\newtheorem{prop}[thm]{Proposition}
\newtheorem{cor}[thm]{Corollary}
\newtheorem{lemma}[thm]{Lemma}
\newtheorem{remark}[thm]{Remark}
\newenvironment{rmk}{\begin{remark}\rm}{\end{remark}}

\numberwithin{equation}{section}

\title[Bundle extensions and a generalised theta divisor]{Geometry of vector bundle extensions and applications to a generalised theta divisor}

\date{}

\author{George H.\ Hitching}
\address{H\o gskolen i Oslo og Akershus, Postboks 4, St.\ Olavs plass, Oslo, Norway} \email{george.hitching@hioa.no}

\thanks{I am grateful to Leibniz Universit\"at Hannover, the University of Oslo, the Korean Institute of Advanced Study and Konkuk University, Seoul, for financial support and hospitality. I acknowledge gratefully the financial support of the Deutsche Forschungsgemeinschaft Schwerpunktprogramm ``Globale Methoden in der Komplexen Geometrie'' and the Norges Forskningsr\aa d ``Strategisk Universitetsprogram i Ren Matematikk''. At the time of submission, the author's affiliation was to H\o gskolen i Vestfold, Bakkenteigen, Norway.}



\begin{document}

\begin{abstract}
Let $E$ and $F$ be vector bundles over a complex projective smooth curve $X$, and suppose that $0 \to E \to W \to F \to 0$ is a nontrivial extension. Let $G \subseteq F$ be a subbundle and $D$ an effective divisor on $X$. We give a criterion for the subsheaf $G(-D) \subset F$ to lift to $W$, in terms of the geometry of a scroll in the extension space $\PP H^{1}(X, \Homom(F, E))$. We use this criterion to describe the tangent cone to the generalised theta divisor on the moduli space of semistable bundles of rank $r$ and slope $g-1$ over $X$, at a stable point. This gives a generalisation of a case of the Riemann--Kempf singularity theorem for line bundles over $X$. In the same vein, we generalise the geometric Riemann--Roch theorem to vector bundles of slope $g-1$ and arbitrary rank.
\end{abstract}

\maketitle

\section{Introduction}

Let $X$ be a complex projective smooth curve of genus $g \geq 2$ and let $E$ and $F$ be vector bundles over $X$. It is well known that isomorphism classes of extensions $0 \to E \to W \to F \to 0$ are parametrised by the cohomology group $H^{1}(X, \Homom(F,E))$, the zero element corresponding to the trivial extension $F \oplus E$. These spaces have been much investigated and used in many contexts. They can be used to cover moduli spaces of vector bundles (see Narasimhan--Ramanan \cite{NR1987}, also \cite{Hit2005b}), giving a useful tool for the analysis thereof (see for example Pauly \cite{Pau2002}). Extension spaces also occur naturally as tangent spaces at smooth points of these moduli spaces, and we will say more about this later. And they have been used in coding theory by Johnsen \cite{Joh2003}, Coles \cite{Col2005} and others.\\
\par
The central point of the present article is the following. Let $W$ be a nontrivial extension of $F$ by $E$, and suppose $\overline{\gamma} \colon G \to F$ is a vector bundle inclusion. Let $D$ be an effective divisor on $X$, and consider the sheaf injection $G(-D) \to G$ induced by the section of $\Ox(D)$ vanishing exactly along $D$. We write $\gamma$ for the composed map $G(-D) \to G \to F$, which is a sheaf injection and a generically injective map of vector bundles. It is often of interest to know when $\gamma$ factorises via a map $G(-D) \to W$. For this, one has:
\begin{lemma} \label{NRkey} The map $\gamma$ factorises via $W$ if and only if the class $\delta(W)$ of the extension belongs to the kernel of the induced map 
\[ \gamma^{*} \colon H^{1}(X, \Homom(F, E)) \to H^{1}(X, \Homom(G(-D), E)). \]
\end{lemma}
\begin{proof} This is a special case of Narasimhan--Ramanan \cite{NR1969}, Lemma 3.1. \end{proof}

Since nontrivial extensions with proportional extension classes are isomorphic as vector bundles, we lose little by working with $\PP H^{1}(X, \Homom(F, E))$ (in order to avoid trivial cases, we will assume that $h^{1}(X, \Homom(F, E)) \geq 1$).
\par
When $E$ and $F$ are line bundles, Lange and Narasimhan in \cite{LNar1983} gave a geometric criterion for the lifting of $\gamma \colon F(-D) \to F$ to an extension $W$, in terms of secants to the natural image of the curve in
\[ |\Kx F E^{-1}|^{*} = \PP H^{0}(X, \Kx \otimes F \otimes E^{*})^{*} \cong \PP H^{1}(X, \Homom(F, E)). \]
Our first aim is to find an analogous criterion for such liftings when $E$ and $F$ may have higher rank. This will allow us to use geometric methods to illuminate questions associated to such liftings. In the second part of the article, we investigate one such question: the geometric description of tangent cones to a generalised theta divisor.\\
\par
Here is a more precise summary of the article. In $\S$\ref{scrolls}, we consider an arbitrary vector bundle $V \to X$ with $h^{1}(X, V) \geq 1$. We describe a map $\psi \colon \PP V \dashrightarrow \PP H^{1}(X, V)$ and recall conditions for it to be an embedding (Theorem \ref{psiemb}). We also discuss some elementary geometry of varieties in projective space.
\par
Next, let $E$, $F$, $G$, $D$ and $\gamma$ be as above. In $\S$\ref{extlift}, we give a geometric criterion (Theorem \ref{liftcrit}) for the lifting of $\gamma$ to an extension $W$ of $F$ by $E$. This is given in terms of the image of the aforementioned map $\psi$ when $V = \Homom(F, E)$. We assume that $\psi$ is an embedding (as it is in the case where we will later apply Theorem \ref{liftcrit}), but at the end of the section we sketch how the criterion can be made sense of when this hypothesis is not satisfied.
\par
In $\S$\ref{gentheta}, we apply Theorem \ref{liftcrit} to the study of a generalised theta divisor. The moduli space ${\mathcal U}(r, r(g-1))$ of semistable bundles of rank $r$ and slope $g-1$ over $X$ has a natural divisor $\Delta$, whose support consists of bundles with nonzero sections. The tangent space to $\U(r, r(g-1))$ at a smooth (equivalently, stable) point $E$ is naturally isomorphic to $H^{1}(X, \Endom(E))$, which parametrises extensions $0 \to E \to \mathbb{E} \to E \to 0$. We use Theorem \ref{liftcrit} to give a description of the projectivised tangent cone to $\Delta$ at $E$ (Theorem \ref{gRKst}), which generalises a case of the Riemann--Kempf singularity theorem (see Griffiths--Harris \cite[Chapter 2]{GH1978}). We use several results of Laszlo \cite{Las1991}.
\par
Furthermore, from Theorem \ref{liftcrit} we deduce a generalisation of the geometric Riemann--Roch theorem (Theorem \ref{ggRRt}), relating the number of sections of a vector bundle $E$ of rank $r$ and slope $g-1$ to the codimension of the linear span of a certain variety in projective space. This result holds even if $E$ is not semistable.
\par
\textbf{Acknowledgements:} I thank Hans-Christian von Bothmer, Insong Choe, Cord Erdenberger, Klaus Hulek, Atanas Iliev, Trygve Johnsen, Peter Newstead, Christian Pauly, Ragni Piene, Kristian Ranestad and Arne B.\ Sletsj{\o}e for enjoyable and valuable discussions. 

\section{Scrolls in projective space} \label{scrolls}

Let $V \to X$ be a vector bundle of rank $r \geq 1$ such that $\Kx \otimes V^*$ has at least one section. We describe a rational map of the scroll $\PP V \to X$ into the projective space $\PP H^{1}(X,V)$. Let $\pi \colon \PP V \to X$ be the projection. We have the following sequence of identifications:
\begin{align*} H^{1}(X,V) &\cong H^{0}(X, \Kx \otimes V^{*} )^{*} \hbox{ by Serre duality} \\
 &\cong H^{0}(X, \Kx \otimes \pi_{*}O_{\PP V}(1) )^{*} \\
 &\cong H^{0}(X, \pi_{*} \left( \pi^{*}\Kx \otimes O_{\PP V}(1) \right) )^{*} \hbox{ by the projection formula} \\
 &\cong H^{0}(\PP V, \pi^{*}\Kx \otimes O_{\PP V}(1) )^{*} \hbox{ by definition of direct image.} \end{align*}
By standard algebraic geometry, we have a map $\psi \colon \PP V \dashrightarrow \PP H^{1}(X, V)$. We write $\Upsilon$ for the line bundle $\pi^{*}\Kx \otimes O_{\PP V}(1) \to \PP V$.
\begin{thm} \label{psiemb}
\begin{enumerate}
\renewcommand{\labelenumi}{(\arabic{enumi})}
\item The map $\psi \colon \PP V \dashrightarrow \PP H^{1}(X, V)$ is an embedding if we have $h^{0}(X, V) = h^{0}(X, V(D))$ for all effective divisors $D$ of degree two on $X$.
\item In particular, $\psi$ is an embedding if $V = \Homom(F, E)$ for semistable bundles $E$ and $F$ with $\mu(F) > \mu(E) + 2$. \end{enumerate} \end{thm}
\begin{proof} Ramanan and Hwang prove (1) in \cite[$\S$3]{HR2004} if $V = \Endom_{0}(E)$, the bundle of trace zero endomorphisms of another bundle $E$, and in fact their argument also applies to arbitrary $V$. As for (2): if $E$ and $F$ are semistable and $\mu(F) > \mu(E) + 2$ then there are no maps $F \to E(D)$ for any divisor $D$ of degree two on $X$. Thus $h^{0}(X, \Homom(F, E)) = 0 = h^{0}(X, \Homom(F, E(D)))$, and we can use (1). \end{proof}

\begin{rmk} Let $L \to X$ be a line bundle. If we identify $X$ with $\PP \Kx L^{-1}$ then, via Serre duality, $\psi$ coincides with the standard map $\phi_{L} \colon X \dashrightarrow |L|^*$. Then it is easy to check that Theorem \ref{psiemb} gives a direct generalisation of the well-known fact that $\phi_L$ is an embedding if and only if $h^{0}(X, L(-x-y)) = h^{0}(X, L) - 2$ for all $x,y \in X$. \end{rmk}

We now recall two facts on varieties in projective space:

\subsection*{Linear span of a subvariety}

Suppose $Y$ is a smooth projective variety and $Z$ a closed subvariety of $Y$, and let $\Upsilon \to Y$ be a line bundle with at least one section. We write $\psi$ for the standard map $Y \dashrightarrow |\Upsilon|^*$. We will describe the linear span of $\psi(Z)$ in $|\Upsilon|^*$.
\par
We write $\cI_Z$ for the ideal sheaf of $Z$. Then $H^{0}(Y, \Upsilon \otimes \cI_{Z})$ is the subspace of global sections of $\Upsilon$ which vanish along $Z$. These correspond to hyperplanes in $|\Upsilon|^*$ which contain the image of $Z$. Let $\Upsilon'$ be the restriction of $\Upsilon$ to $Z$. There is a natural exact sequence
\[ 0 \to H^{0}(Y, \Upsilon \otimes \cI_{Z} ) \to H^{0}(Y, \Upsilon) \xrightarrow{q} H^{0}(Z, \Upsilon') \to \cdots
\]\begin{prop} The linear span of $\psi(Z)$ in $|\Upsilon|^*$ coincides with the projectivised kernel of $H^{0}(Y, \Upsilon)^{*} \to H^{0}(Y, \Upsilon \otimes \cI_{Z})^*$. \label{linspan} \end{prop}
\begin{proof} Write $\Pi$ for the image of $q$ in $H^{0}(Z, \Upsilon')$; then
\[ \Pi^{*} \cong \Ker \left( H^{0}(Y, \Upsilon)^{*} \to H^{0}(Y, \Upsilon \otimes \cI_{Z})^{*} \right). \]
We write $\psi' \colon Z \dashrightarrow | \Upsilon' |^*$ for the map to projective space determined by $\Upsilon'$. It is easy to check that there is a commutative diagram
\[ \xymatrix{ & Z \ar @^{(->}[r] \ar @{-->}[dl]_{\psi'} & Y \ar @{-->}[dd]^{\psi} \\
|\Upsilon'|^{*} \ar @_{-->}[dr] \ar @_{-->}[drr]^{\PP q^*} & & \\
 & \PP \Pi^{*} \ar @^{(->}[r] & |\Upsilon|^{*} } \]
and this shows that the linear span of $\psi(Z)$ is contained in $\PP \Pi^*$. Moreover, since $\psi'(Z)$ is nondegenerate in $| \Upsilon' |^*$ and the restriction $H^{0}(Z, \Upsilon')^{*} \to \Pi^*$ is surjective, the image of $Z$ in $\PP \Pi^*$ is also nondegenerate. This proves the proposition. \end{proof}

\subsection*{Osculating spaces} \label{osc}

References for this subject include Piene \cite{Pie1976}, Piene--Tai \cite{PT1990} and Lanteri et al \cite{LMP2006}. Let $Y$ be a smooth projective variety and $\Upsilon \to Y$ a very ample line bundle. We denote $\psi \colon Y \hookrightarrow |\Upsilon|^*$ the map to projective space defined by $\Upsilon$. For $y \in Y$ and $k \geq 0$, the \textsl{$k$th osculating space to $Y$ at $y$} is defined as the projective linear subspace of $|\Upsilon|^*$ spanned by the forms on $H^{0}(Y,\Upsilon)^*$ defined by differential operators of order at most $k$ at $y$. We denote it $\Osc^{k}(Y, y)$. For large $k$, it will fill up all of $|\Upsilon|^*$.
\par
By choosing a system of local coordinates near $y$, we see that those sections of $\Upsilon$ which are annihilated by all differential operators of order at most $k$ are exactly those which vanish to order at least $k+1$ at $y$; precisely, those whose images in the stalk $\Upsilon_y$ belong to $( \Upsilon \otimes \cI_{y}^{k+1} )_y$, where $\cI_y$ is the ideal sheaf of $y$. Thus we have:
\begin{prop} \label{oscsp} The $k$th osculating space to $Y$ at $p$ coincides with
\[ \Osc^{k}(Y,y) = \PP \Ker \left( H^{0}(Y, \Upsilon)^{*} \to H^{0}(Y, \Upsilon \otimes \cI_{y}^{k+1})^{*} \right). \]
\end{prop}
This description of $\Osc^{k}(Y,y)$ will be useful in what follows.

\section{Extensions, lifting and geometry} \label{extlift}

Let $E$ and $F$ be vector bundles over $X$. Throughout this section, we will suppose that the map $\psi \colon \PP \Homom(F, E) \dashrightarrow \PP H^{1}(X, \Homom(F, E))$ defined in the last section is an \emph{embedding}. At the end of the section, we will briefly discuss what happens more generally.\\
\par
As before, let $\gamma$ be the generically injective vector bundle map defined by the composition
\[ G(-D) \to G \xrightarrow{\ogamma} F, \]
where $\ogamma$ is a vector bundle inclusion and $G(-D) \to G$ is induced by the section of $\Ox(D)$ vanishing precisely along $D$.
\par
We now define some loci in $\PP H^{1}(X, \Homom(F, E))$. Since $\ogamma$ is a vector bundle injection, the kernel of the induced map
\[ \ogamma^{*} \colon \Homom(F, E) \to \Homom(G, E) \]
is a vector subbundle $\uKer(\ogamma^{*})$ of $\Homom(F, E)$, of rank $(\rank(F) - \rank(G)) \cdot \rank(E)$. If we write $H$ for $\uCoker( \ogamma )$, a vector bundle, then $\uKer ( \ogamma^{*} ) \cong \Homom( H, E )$. We obtain a subscroll $\PP \uKer( \ogamma^{*})$ of $\PP \Homom(F, E)$, which is nonempty if and only if $\rank(F) > \rank(G)$.
\par
Next, we write $D = \sum_{i=1}^{n} k_{i} x_i$, where the $x_i$ are distinct and each $k_{i} \geq 1$. We define $N(\gamma)$ to be the union of the following loci in $\PP H^{1}(X, \Homom(F, E))$:
\begin{itemize}
\item the subscroll $\PP \uKer(\ogamma^{*})$ of $\PP \Homom(F, E)$, and
\item the union over all $i = 1, \ldots , n$ of
\[ \bigcup_{\nu \in \PP \Homom(F, E)|_{x_i}} \Osc^{k_{i}-1} \left( \PP \Homom(F, E), \nu \right). \]
\end{itemize}
We can now state the main result of this section, which generalises the idea of Lange--Narasimhan \cite[Proposition 1.1]{LNar1983}.

\begin{thm} Let $W$ be a nontrivial extension of $F$ by $E$. Then the map $\gamma \colon G(-D) \to F$ lifts to $W$ if and only if $\langle \delta(W) \rangle$ belongs to the linear span of the locus $N(\gamma)$ in $\PP H^{1}(X, \Homom(F, E))$. \label{liftcrit} \end{thm}
The proof of this theorem will occupy the remainder of this section. Let us briefly indicate the principle before getting into the details: By Lemma \ref{NRkey}, we need to show that the linear span of $N(\gamma)$ is equal to 
\[ \PP \Ker \left( \gamma^{*} \colon H^{1}(X, \Homom(F, E)) \to H^{1}(X, \Homom(G(-D), E)) \right) \]
or equivalently, by Serre duality,
\[ \PP \Ker \left( \gamma^{*} \colon H^{0}(X, \Kx \otimes F \otimes E^{*} )^{*} \to H^{0}(X, \Kx \otimes G(-D) \otimes E^{*} )^{*} \right). \]
In the last section we saw that $H^{0}(X, \Kx \otimes F \otimes E^{*})$ is identified with $H^{0}(\PP \Homom(F, E), \Upsilon )^*$, where $\Upsilon$ is a certain line bundle. We will show that under this identification, $H^{0}(X, \Kx \otimes G(-D) \otimes E^{*} )$ corresponds to a subspace of $H^{0}(\PP \Homom(F, E), \Upsilon)$ of sections vanishing along the various components of $N(\gamma)$ with appropriate multiplicities. Then we use Propositions \ref{linspan} and \ref{oscsp} to conclude.\\
\par
We begin by assembling some technical results. Let $V \to X$ be any vector bundle. We will study the connection between sections of $\Kx \otimes V^{*} \to X$ and those of $\Upsilon \to \PP V$ in more detail. Let $Q$ be a subbundle of $V$. We write $Q^{\perp}$ for the orthogonal complement of $Q$ in $V^*$, which is defined by the exact sequence $0 \to Q^{\perp} \to V^{*} \to Q^{*} \to 0$.
\begin{lemma} \label{sectcorrA} Via the identification
\begin{equation} H^{0}(X, \Kx \otimes V^{*} ) \xrightarrow{\sim} H^{0}(\PP V, \Upsilon) \label{sectcorr} \end{equation}
described in $\S$\ref{scrolls}, the subspace $H^{0}(X, \Kx \otimes Q^{\perp})$ of $H^{0}(X, \Kx \otimes V^{*})$ corresponds to the space of global sections of $\Upsilon$ vanishing along $\PP Q$. \end{lemma}
\begin{proof} For any $x \in X$, a section $s$ of $\Kx \otimes V^*$ restricts to a linear form on $V|_x$ with values in the line $\Kx|_x$. To evaluate the corresponding section $\tilde{s}$ of $\Upsilon \to \PP V$ at a point $\nu \in \PP V|_x$, we restrict $s$ to the line $\nu \subseteq V|_x$ and obtain an element of $\Homom( \nu, \Kx|_{x}) = \C$. Thus $\tilde{s}$ vanishes at $\nu$ if and only if
\[ s(x) \in \Ker \left( \Kx \otimes V^{*}|_{x} \to \Kx|_{x} \otimes \nu^{*} \right) . \]
In the same way, we see that $\tilde{s}$ vanishes at all points of $\PP Q|_x$ if and only if
\[ s(x) \in \Ker \left(\Kx \otimes V^{*} \to \Kx \otimes Q^{*} \right)|_{x}, \]
that is, $s(x) \in \Kx \otimes Q^{\perp}|_x$. Therefore, sections of $\Upsilon \to \PP V$ vanishing along the whole of $\PP Q$ correspond to sections of $\Kx \otimes V^{*} \to X$ which are everywhere $\Kx \otimes Q^{\perp}$-valued. \end{proof}
We adapt this lemma to the situation in which we will need it:
\begin{cor} Suppose that $V = \Homom(F, E)$ for vector bundles $E$ and $F$, and that $\ogamma \colon G \hookrightarrow F$ is a vector bundle inclusion. We consider again the subbundle
\[ \uKer \left( \ogamma^{*} \colon \Homom(F, E) \to \Homom (G, E) \right) \]
of $\Homom(F, E)$. Then sections of $\Upsilon \to \PP \Homom(F, E)$ vanishing along the subscroll $\PP \uKer(\ogamma^{*})$ correspond to sections of $\Kx \otimes F \otimes E^{*} \to X$ with values in $\Kx \otimes G \otimes E^*$ at all points. \label{sectcorrcor} \end{cor}
\begin{proof} In view of Lemma \ref{sectcorrA}, it suffices to show that the orthogonal complement of $\uKer(\ogamma^{*})$ in $F \otimes E^*$ is $G \otimes E^*$. Recall that $\uKer ( \ogamma^{*} ) \cong H^{*} \otimes E$ where $H = \uCoker ( \ogamma )$. Thus we have an exact sequence
\[ 0 \to \uKer( \ogamma^{*} ) \to F^{*} \otimes E \to G^{*} \otimes E \to 0. \]
Dualising, we see that $\uKer ( \ogamma^{*} )^{\perp} = G \otimes E^*$, as required. \end{proof}

Next, let $x$ be a point of $X$. We denote $m$ the maximal ideal of $\Oxx$, and for any $\nu \in \PP V$ we write $\cI_{\nu}$ for the maximal ideal of ${\mathcal O}_{\PP V, \nu}$.
\begin{lemma} For each $k \geq 0$, sections of $\Kx \otimes V^*$ belonging to $m^{k} \left( \Kx \otimes V^{*} \right)_x$ correspond via (\ref{sectcorr}) to sections of $\Upsilon$ belonging to $\cI_{\nu}^k$ for \emph{all} $\nu \in \PP V|_x$. \label{sectcorrB} \end{lemma}
\begin{proof} Let $s$ be a section of $\Kx \otimes V^*$, and let $j$ be the unique integer such that
\[ s_{x} \in m^{j} \left( \Kx \otimes V^{*} \right)_{x} \backslash m^{j+1} \left( \Kx \otimes V^{*} \right)_{x}. \]
Then for any uniformiser $z$ at $x$, the restriction of $s$ to a suitable neighbourhood $U$ of $x$ has the form $z^{j} \cdot t$ for some section $t$ of $\Kx \otimes V^*$ with \nolinebreak
\[ t_{x} \in (\Kx \otimes V)_{x} \backslash m(\Kx \otimes V)_{x}. \]
The section $t$ is well-defined up to a unit in $O_{X,x}$.
\par
We describe the corresponding section $\tilde{s}$ of $\Upsilon$ more precisely. Firstly, we consider the section $\pi^{*}s$ of $\pi^{*} \left( \Kx \otimes V^{*} \right)$. By Hartshorne \cite[II.7.11(b)]{Har1977}, there is a canonical surjection $\pi^{*} V^{*} \to O_{\PP V}(1)$, which is identified on each fibre with restriction of a global section of $O_{\PP^{r-1}}(1)$ to each point of $\PP^{r-1}$ in turn. The image of $\pi^{*}s$ under the induced map
\[ \pi^{*} \left( \Kx \otimes V^{*} \right) \to \pi^{*} \Kx \otimes O_{\PP V}(1) = \Upsilon \]
is the desired section $\tilde{s}$ of $\Upsilon$.
\par
Suppose now that $j \geq k$. Now $\pi^{*} z$ is a local function on $\PP V$, which belongs to $\cI_{\nu} \backslash \cI_{\nu}^2$ for all points $\nu$ of $\PP V|_x$. The section $\tilde{s}$ of $\Upsilon$ determined as above by $s$ is a multiple of $(\pi^{*} z)^k$, and therefore $\tilde{s} \in \cI_{\nu}^{k}\Upsilon_{\nu}$ for all $\nu \in \PP V|_x$.
\par
Conversely, suppose that $j < k$. Choose any $\nu_{0} \in \PP V|_x$ such that $t|_{\nu_0}$ is nonzero. Then $\tilde{s}$ lies outside $\cI_{\nu_0}^{j+1} \Upsilon_{\nu_0}$, and in particular outside $\cI_{\nu_0}^{k} \Upsilon_{\nu_0}$. \end{proof}

The last technical tool we need is a result in linear algebra. Let $N$ be a finite-dimensional vector space and $\{ N_{\lambda} : \lambda \in I \}$ a collection of subspaces of $N$ indexed by a set $I$. Write $N_0$ for the intersection of all the $N_{\lambda}$.
\begin{prop} The kernel of the restriction map $N^{*} \to N_{0}^*$ coincides with the linear span $S$ of the kernels of the restriction maps $N^{*} \to N_{\lambda}^*$ for all $\lambda \in I$. \label{kernels} \end{prop}
\begin{proof} Firstly, dualising the exact sequence $0 \to N_{\lambda}^{\perp} \to N^{*} \to N_{\lambda}^{*} \to 0$, we see that $\left( N_{\lambda}^{\perp} \right)^{\perp}$ coincides with $N_{\lambda}$ under the canonical identification $N = (N^{*})^*$. For $v \in N$, we have
\begin{align*} \hbox{$\phi(v) = 0$ for all $\phi \in S$} \quad \Longleftrightarrow & \quad \phi(v) = 0 \hbox{ for all $\phi \in N_{\lambda}^{\perp}$ and all $\lambda$} \\
 \Longleftrightarrow & \quad v \in \left( N_{\lambda}^{\perp} \right)^{\perp} = N_{\lambda} \hbox{ for all $\lambda$} \\
 \Longleftrightarrow & \quad v \in N_{0}. \end{align*}
Thus $S = \Ker( N^{*} \to N^{*}_{0})$, as required. \end{proof}
Now we have all the ingredients for our main result on liftings:\\
\\
\textit{Proof of Theorem \ref{liftcrit}.} As mentioned before, by Lemma \ref{NRkey}, we need to show that the linear span of $N(\gamma)$ is equal to
\[ \PP \Ker \left( \gamma^{*} \colon H^{1}(X, \Homom(F, E) ) \to H^{1}(X, \Homom(G(-D), E)) \right) \]
or, equivalently, by Serre duality,
\[ \PP \Ker \left( \gamma^{*} \colon H^{0}(X, \Kx \otimes F \otimes E^{*})^{*} \to H^{0}(X, \Kx \otimes G(-D) \otimes E^{*})^{*} \right). \]
Now the image of $H^{0}(X, \Kx \otimes G(-D) \otimes E^{*})$ in $H^{0}(X, \Kx \otimes F \otimes E^{*} )$ is equal to the intersection
\[ \ogamma \left( H^{0}(X, \Kx \otimes G \otimes E^{*}) \right) \cap \left( \bigcap_{i=1}^{n} H^{0} \left( X, \Kx \otimes F( - k_{i}x_{i}) \otimes E^{*} \right) \right). \]
By Proposition \ref{kernels}, therefore, $\Ker( \gamma^{*} )$ is the linear span of the union of
\begin{equation} \Ker\left( H^{0} ( X, \Kx \otimes F \otimes E^{*} )^{*} \to H^{0} ( X, \Kx \otimes G \otimes E^{*})^{*} \right) \label{AA} \end{equation}
together with the union of
\begin{equation} \Ker \left( H^{0}(X, \Kx \otimes F \otimes E^{*} )^{*} \to H^{0} \left( X, \Kx \otimes F \left( - k_{i}x_{i} \right) \otimes E^{*} \right)^{*} \right) \label{BB} \end{equation}
over all $i = 1, \ldots, n$.
\par
By Corollary \ref{sectcorrcor}, the space $H^{0}(X, \Kx \otimes G \otimes E^{*})$ corresponds to the space of global sections of $\Upsilon \to \PP \Homom(F, E)$ which vanish along the subscroll $\PP \uKer (\ogamma^{*})$. Therefore the space (\ref{AA}) is identified with
\begin{equation} \Ker \left( H^{0}( \PP \Homom(F, E), \Upsilon )^{*} \to H^{0}( \PP \Homom(F, E), \Upsilon \otimes \cI_{\PP \uKer( \ogamma^{*} )} )^{*} \right). \label{sectA} \end{equation}
Next, by Lemma \ref{sectcorrB}, the space $H^{0} \left( X, \Kx \otimes F ( - k_{i}x_{i}) \otimes E^{*} \right)$ corresponds to
\[ \bigcap_{ \nu \in \PP \Homom(F, E)|_{x_i}} H^{0} \left( \PP \Homom(F, E), \Upsilon \otimes \cI_{\nu}^{k_i} \right). \]
Thus by Proposition \ref{kernels}, the space (\ref{BB}) corresponds to the linear span of
\begin{equation} \bigcup_{\nu \in \PP \Homom(F, E)|_{x_i}} \Ker \left( H^{0}( \PP \Homom(F, E), \Upsilon )^{*} \to H^{0}( \PP \Homom(F, E), \Upsilon \otimes \cI_{\nu}^{k_i})^{*} \right) \label{sectB} \end{equation}

Putting all this together, $\Ker( \gamma^{*} )$ is exactly the linear span of (\ref{sectA}) and the union of the loci (\ref{sectB}) for $i = 1, \ldots, n$.
\par
Now we projectivise: (\ref{sectA}) becomes $\Span \left( \PP \uKer (\ogamma^{*}) \right)$ by Proposition \ref{linspan} and (\ref{sectB}) becomes
\[ \bigcup_{\nu \in \PP \Homom(F, E)|_{x_i}} \Osc^{k_{i} - 1} \left( \PP \Homom(F, E), \nu \right) \]
by Proposition \ref{oscsp}. Thus $\PP \Ker(\gamma^{*})$ is exactly the linear span of $N(\gamma)$.\\
\\
This completes the proof of Theorem \ref{liftcrit}. \qed

\subsection*{If $\psi$ is not an embedding}

If we drop the assumption that $\psi$ be an embedding, then what we have shown in Theorem \ref{liftcrit} is the following: \textit{The sheaf injection $\gamma \colon G(-D) \to F$ factorises via $W$ if and only if $\delta(W)$ belongs to the linear span of the union of the following:}
\begin{itemize}
\item $\Ker \left( H^{0} ( \PP \Homom(F, E), \Upsilon)^{*} \to H^{0}(\PP \Homom(F, E), \Upsilon \otimes \cI_{\PP \uKer( \ogamma^{*} )})^{*} \right)$, and
\item the union over $i = 1, \ldots, n$ and $\nu \in \PP \Homom(F, E)|_{x_i}$ of
\[ \Ker \left( H^{0} ( \PP \Homom(F, E), \Upsilon)^{*} \to H^{0}(\PP \Homom(F, E), \Upsilon \otimes \cI_{\nu}^{k_i})^{*} \right). \]
\end{itemize}
If $\psi$ does not fail too badly to be an embedding (for example, if it contracts some pairs of points, if the differential fails to be injective at some points, or even if $\psi$ fails to be defined at some points) then the above spaces can be interpreted geometrically in an obvious way using Lemmas \ref{linspan} and \ref{oscsp}. Notice for example that in the statement of Lemma \ref{linspan}, we did not assume that $\psi$ was an embedding. We will give one example which shows the effect such phenomena have on the behaviour of the extensions. Suppose $X$ is a hyperelliptic curve and we are considering rank two extensions of the form
\[ 0 \to \Ox \to W \to \Ox \to 0 \]
with classes in $H^{1}(X, \Ox)$. Let $D = \sum_{i=1}^{n} k_{i}x_i$ be an effective divisor on $X$, and let $\gamma$ be the inclusion $\Ox(-D) \to \Ox$. Here $\ogamma$ is the identity $\Ox \xrightarrow{\sim} \Ox$, so the scroll $\PP \uKer(\ogamma^{*})$ is empty, and $N(\gamma)$ consists of
\[ \bigcup_{i=1}^{n} \Osc^{k_{i}-1}(\phi(X), \phi(x_{i})). \]
Thus the span of $N(\gamma)$ is precisely the secant space to $X$ spanned by $D$. By Theorem \ref{liftcrit}, the map $\gamma$ factorises via $W$ if and only if $\delta(W)$ lies on this secant. (This also follows from Lange--Narasimhan \cite[Proposition 1.1]{LNar1983}.)
\par
Now write $\gamma' \colon \Ox(-\iota(D)) \to \Ox$ where $\iota \colon X \to X$ is the hyperelliptic involution. The locus $N(\gamma')$ is the secant to $\psi(X)$ spanned by the image of the divisor $\iota(D)$. But this is precisely the image of $D$. Therefore we have
\[ N(\gamma) = N(\gamma'). \]
Hence, by Theorem \ref{liftcrit}, any extension $W$ to which $\gamma$ lifts must also admit a lifting from $\gamma'$, and vice versa.
\par
Thus, the fact that the map $\PP \Homom(\Ox, \Ox) \to \PP H^{1}(X, \Ox)$ fails to be injective is reflected in a natural way in the properties of the extensions. (This example should be compared with Lange--Narasimhan \cite[p.\ 59]{LNar1983}, especially if $\deg(D) = 1$.)

\section{Tangent cones of a generalised theta divisor} \label{gentheta}

{\footnotesize \textbf{Note:} To readers following the reference in Lemma 4.5 of the article \cite{HJ2008} of Trygve Johnsen and the present author: Please note that the reference is to a early version of the present article, which can be found online at \texttt{http://arxiv.org/abs/math/0610970v3} . I \nolinebreak apologise for this inconvenience.}\\
\par
In this section we use Theorem \ref{liftcrit} to generalise two well-known results on line bundles over curves to bundles of higher rank: the Riemann--Kempf singularity theorem and the geometric Riemann--Roch theorem. We begin with a brief review of these results in the line bundle case.\\
\par
Let $X$ be a curve of genus $g \geq 3$. We assume for simplicity that $X$ is nonhyperelliptic, so that the canonical map $\phi \colon X \hookrightarrow |\Kx|^*$ is an embedding. We denote $J^{g-1}$ the Jacobian variety parametrising line bundles of degree $g-1$ over $X$, and we write $W_{g-1}$ for the natural divisor on $J^{g-1}$ whose support consists of bundles with sections.
\par
Let $L \to X$ be a line bundle of degree $g-1$ with $h^{0}(X,L) \geq 1$. Recall the \emph{geometric Riemann--Roch theorem}:
\begin{thm} \label{gRRt} The number $h^{0}(X, L)$ is equal to $\codim ( \Span(\phi(D)), |\Kx|^{*} )$.\end{thm}
\begin{proof} See Griffiths--Harris \cite{GH1978}, chapter 2. \end{proof}

Next, we write $h^{0}(X, L) =: n$. By the Riemann singularity theorem, $L$ is a point of multiplicity $n$ in $W_{g-1}$. Recall that the projectivised tangent space to $J^{g-1}$ at any point is isomorphic to $|\Kx|^*$. The \emph{Riemann--Kempf singularity theorem} gives us a geometric description of the tangent cones of $W_{g-1}$:

\begin{thm} \label{RKst} The projectivised tangent cone to $W_{g-1}$ at $L$ is the union of the projective $(g-1-n)$-planes spanned by the images by $\phi$ of all the effective divisors in $|L|$. \end{thm}
\begin{proof} See Griffiths--Harris \cite{GH1978}, chapter 2. \end{proof}

\subsection{A generalisation of the Riemann--Kempf singularity theorem}

We begin by describing the objects which will replace $J^{g-1}$ and $W_{g-1}$. We recall some facts from Laszlo \cite{Las1991}:\\
\par
Let $\U := {\mathcal U}(r, r(g-1))$ be the moduli space of semistable vector bundles of rank $r$ and degree $r(g-1)$. This is a projective irreducible normal variety of dimension $r^{2}(g-1)+1$. It has a distinguished divisor $\Delta$, whose support consists of bundles admitting at least one independent section (which clearly coincides with $W_{g-1}$ if $r = 1$). It is well known that the tangent space $T_{E} \U$ to $\U$ at a smooth point $E$ is isomorphic to $H^{1}(X, \Endom(E))$.\\
\\
We now state an important fact:

\begin{thm} \label{psigenemb} Suppose that $X$ is nonhyperelliptic and of genus $g \geq 5$. Then the map $\psi \colon \PP \Endom(E) \dashrightarrow \PP H^{1}(X, \Endom(E)) = \PP T_{\U}|_E$ defined in $\S$\ref{scrolls} is an embedding for generic $E \in \Delta$. \end{thm}
\begin{proof} This theorem is proven by Hwang and Ramanan in \cite[$\S$3]{HR2004} if $\Endom(E)$ is replaced with $\Endom_{0}(E)$, the bundle of trace zero endomorphisms of $E$, and the result we require will follow easily from their work. By Proposition \ref{psiemb}, we need to show that for generic $E \in \Delta$ we have
\[ h^{0}(X, \Endom(E)) = h^{0}(X, \Homom(E, E(D))) \]
for all effective divisors $D$ of degree two over $X$. Recall the natural direct sum decomposition
\[ \Endom(E) = \Ox \oplus \Endom_{0}(E). \]
Since $X$ is nonhyperelliptic, we have $h^{0}(X, \Ox) = 1 = h^{0}(X, \Ox(D))$ for all $D \in \Sym^{2}X$. Furthermore,
\[ h^{0}(X, (\Endom_{0}(E))(D)) = 0 = h^{0}(X, \Endom_{0}(E)) \]
for generic $E \in \U$ by \cite[Proposition 3.2]{HR2004}. Therefore, for generic $E \in \U$, the map $\psi$ is indeed an embedding.
\par
We must show that in fact the statement is true for a general bundle in the divisor $\Delta$. To see this, note that a twist $F \otimes L$ of a general bundle $F \in \U$ by a suitable line bundle $L$ of degree zero will belong to $\Delta$ (precisely: an $L$ belonging to the \emph{theta divisor} of $F$. See for example Laszlo \cite[$\S$I.2]{Las1991}). But $\Endom(F \otimes L) \cong \Endom(F)$ for any line bundle $L$. The result follows. \end{proof}

Let $E \in \U$ be a stable bundle with $h^{0}(X, E) = n \geq 1$ and such that $\psi$ is an embedding. As in Laszlo \cite{Las1991} (see also Narasimhan--Ramanan \cite{NR1975} and Pauly \cite{Pau2003}), we can find an \'etale affine neighbourhood $S = \Spec (A)$ of $E$ in $\U$ and a family of stable vector bundles $\mathbf{E}$ over $S \times X$ such that for each $F \in S$, we have $\mathbf{E}|_{\{ F \} \times X} \cong F$, together with a homomorphism $\mu \colon M \to N$ of flat $A$-modules of finite type such that for all $A$-modules $P$, by functoriality,
\begin{multline*} H^{0}(S \times X, \mathbf{E} \otimes_{A} P) \cong \Ker \left( \mu \otimes \Iden_{P} \right) \quad \hbox{and} \\ H^{1}(S \times X, \mathbf{E} \otimes_{A} P) \cong \Coker \left( \mu \otimes \Iden_{P} \right). \end{multline*}
Moreover, shrinking $S$, we can suppose that $M$ and $N$ are free $A$-modules and $\mu|_E$ is the zero homomorphism. The divisor $\Delta$ is given on $S$ as $\left( \det ( \mu ) \right)$. Laszlo has given the following generalisation of the Riemann singularity theorem:
\begin{thm} \label{Lst} The multiplicity of $\Delta$ at $E$ is equal to $h^{0}(X, E)$. \end{thm}
\begin{proof} See \cite[Th\'eor\`eme II.10]{Las1991}. \end{proof}

Now let $s$ be a nonzero section of $E$. We regard $s$ as a vector bundle map $\Ox \to E$. As such, it factorises $\Ox \to \Ox(D) \xrightarrow{\overline{s}} E$, where
\[ D = \sum_{i=1}^{n}k_{i}x_i \]
is the divisor of zeroes of $s$, and $\Ox(D)$ is the vector subbundle of $E$ generated by $s(\Ox)$. (If $\rank (E) \geq 2$ then we expect $D$ to be zero for general $s$.) The saturated map $\overline{s} \colon \Ox(D) \to E$ is a vector bundle injection, so
\[ \uKer \left( \overline{s}^{*} \colon \Endom(E) \to \Homom(\Ox(D), E) \right) \]
is a vector subbundle of $\Endom(E)$, of rank $r(r - 1)$. Setting $F = E$ and $G = \Ox(D)$, and $\gamma = s$, we are in the situation of Theorem \ref{liftcrit}. Thus we can define the locus $N(s)$ in $\PP H^{1}(X, \Endom(E))$ as in $\S$\ref{extlift}, as the union of the following subvarieties of $\PP H^{1}(X, \Endom(E))$:
\begin{itemize}
\item the subscroll $\PP \uKer(\overline{s}^{*})$ of $\PP \Endom(E)$, and
\item the union over $i = 1, \ldots , n$ and $\nu \in \PP \Endom(E)|_{x_i}$ of the osculating spaces $\Osc^{k_{i}-1} \left( \PP \Endom(E) , \nu \right)$. \end{itemize}

\begin{thm}[Generalised Riemann--Kempf singularity] \label{gRKst} The projectivised tangent cone to $\Delta$ at $E$ is, set-theoretically, the union over all nonzero $s \in H^{0}(X, E)$ of the linear spans of $N(s) \subset \PP T_{E} \U$. \end{thm}
\begin{proof} Let $A$, $S=\Spec(A)$, $M$ and $N$ be as above and write $\m$ for the maximal ideal of the point $E$ in $A$. By Narasimhan--Seshadri \cite[Lemma 2.1 (ii)]{NS1965}, near $E$ the variety $\U$ looks like an analytic open set in $H^{1}(X, \Endom(E))$. Thus we have flat structures on $\U$ at $E$ to all orders (see Kempf \cite{Kem1986}).
\par
By Theorem \ref{Lst}, the function $\det(\mu)$ belongs to $\m^{n} \backslash \m^{n+1}$. Therefore, the tangent cone in which we are interested is defined by $\det(\mu) \modulo \m^{n+1}$, which we regard as a function $f_n$ on $T_{E} \U$ via the flat structure. We then notice, following Laszlo \cite[$\S$II]{Las1991}, that (for a compatible choice of flat structures of orders $1$ and $n$) we have
\[ \det (\mu) \modulo \m^{n+1} = \det \left( \mu \modulo \m^{2} \right) \modulo \m^{n+1}. \]
As in \cite{Las1991}, we will interpret $\mu \modulo \m^2$ in terms of cup products.
\par
Since $\mu|_E$ is zero, $\mu \modulo \m^2$ is a matrix of elements of $\m / \m^2$. Therefore, we can contract it with an element $v$ of
\[ \left( \m / \m^{2} \right)^{*} = T_{E} S \cong T_{E} \U \cong H^{1}(X, \Endom(E)) \]
to obtain a matrix of scalars $(\mu \: \modulo \: \m^{2})(v)$. By \cite[Lemme II.5]{Las1991}, this matrix can be identified with the cup product map $\cdot \cup v \colon H^{0}(X, E) \to H^{1}(X, E)$. Therefore $f_{n}(v) = 0$ if and only if
\[ \det \left( \cdot \cup v \colon H^{0}(X, E) \to H^{1}(X, E) \right) \]
is zero. Thus $\PP \Cone(\Delta, E)$ coincides (set-theoretically) with
\[ \PP \{ v \in H^{1}(X, \Endom(E)) : \det ( \cdot \cup v ) = 0 \}, \]
in other words, the image of the set of $v$ such that the cup product map $\cdot \cup v$ has a kernel.
\par
Now we use the link with extensions. Classes $v \in H^{1}(X, \Endom(E))$ para-metrise extensions $0 \to E \to \mathbb{E}_{v} \to E \to 0$. The coboundary map in the cohomology sequence
\[ 0 \to H^{0}(X, E) \to H^{0}(X, \mathbb{E}_{v}) \to H^{0}(X, E) \to H^{1}(X, E) \to \cdots \]
is none other than cup product by $v$. But this shows that $\cdot \cup v$ has a kernel if and only if a nonzero section of $E$ lifts to the extension $\mathbb{E}_v$. Thus, we have another set-theoretic description of the tangent cone as \emph{the set of those $v$ defining extensions $\mathbb{E}_v$ to which at least one nonzero section lifts from the quotient copy of $E$}. Hence the projectivised tangent cone to $\Delta$ at $E$ is
\[ \PP \left( \bigcup_{s \in H^{0}(X, E)} \left\{ v \in H^{1}(X, \Endom(E)) : \hbox{ $s$ lifts to the extension $\mathbb{E}_{v}$} \right\} \right). \]
But by Theorem \ref{liftcrit}, for each $s \in H^{0}(X, E)$, the locus
\[ \PP \left\{ v \in H^{1}(X, \Endom(E)) : \hbox{ $s$ lifts to the extension $\mathbb{E}_{v}$} \right\} \]
coincides with the linear span of $N(s)$. Thus $\PP \Cone ( \Delta, E)$ is the union of the linear spans of all the $N(s)$, and we are done. \end{proof} 

\begin{rmk} Our description of the tangent cone in terms of cup products can also be explained in terms of deformations of $E$ (compare with Mukai \cite[proof of Proposition 2.6]{Muk1995}). Recall that a tangent vector to $S$ at $E$ is a morphism $\Spec(\C[\varepsilon]) \to S$ which sends the closed point $(\varepsilon)$ to $E$. Such a morphism is determined by a ring homomorphism $\nu \colon A \to \C[\varepsilon]$ satisfying $\nu^{-1}(\varepsilon) = \m$. We then obtain a deformation of $E$ by pulling back the family $\mathbf{E}$ to $X \times \Spec(\C[\varepsilon])$. In particular, we get a short exact sequence
\[ 0 \to E \otimes_{A}(\varepsilon) \to E \otimes_{A} \C[\varepsilon] \to E \otimes_{A} \C \to 0, \]
which naturally yields an extension of $E$ by $E$. With this interpretation, one expects the tangent cone to $\Delta$ at $E$ to correspond to deformations of $E$ which preserve some section of $E$; that is, having nonzero sections lifting from the quotient copy of $E$. But these are exactly those $v$ such that cup product by $v$ has a kernel. \end{rmk}

\subsection{Geometric Riemann--Roch for bundles of higher rank} \label{grr}

Here we show that as in the line bundle case, the number of independent sections of the bundle $E$ is measured by the degeneration in the linear span of a subvariety of $\PP H^{1}(X, \Endom(E))$.

\begin{thm}[Generalised geometric Riemann--Roch] \label{ggRRt} Let $E \to X$ be a vector bundle of rank $r$ and degree $r(g-1)$ with $h^{0}(X, E) \geq 1$. Then
\[ h^{0}(X, E) = \codim \left( \Span ( N(s) ), \PP H^{1}(X, \Endom(E)) \right) \]
for any nonzero section $s$ of $E$. \end{thm}
\begin{proof} For any nonzero section $s \colon \Ox \to E$, the induced map
\[ s^{*} \colon H^{1}(X, \Endom(E)) \to H^{1}(X, E) \]
is surjective, for example since it is Serre dual to the inclusion
\[ H^{0}(X, \Kx \otimes \Ox \otimes E^{*}) \hookrightarrow H^{0}(X, \Kx \otimes E \otimes E^{*}). \]
Therefore $\dim (\Ker(s^{*})) = h^{1}(X, \Endom(E)) - h^{1}(X,E)$. Thus
\[ \dim (\Span (N(s))) = h^{1}(X, \Endom(E)) - h^{1}(X,E) - 1 \]
by Theorem \ref{liftcrit}. Since $\chi(E) = 0$, we have $h^{1}(X,E) = h^{0}(X,E)$, so
\begin{align*} h^{0}(X, E) &= \left( h^{1}(X, \Endom(E)) - 1 \right) - \dim \Span (N(s)) \\
 & = \codim \left( \Span (N(s)), \PP H^{1}(X, \Endom(E)) \right), \end{align*}
as required. \end{proof}
\begin{rmk} This theorem holds even if $E$ is not semistable. Of course, it may then happen that $h^{1}(X, \Endom(E)) > r^{2}(g-1) + 1$. \end{rmk}

\subsection{The line bundle case}

We conclude by examining the statements of Theorems \ref{gRKst} and \ref{ggRRt} when $E$ is a line bundle.
\par
Suppose $X$ is nonhyperelliptic and let $L \to X$ be a line bundle of degree $g-1$. Then the scroll $\PP \Endom(L)$ is simply $\PP \Ox = X$, and $\psi$ is the canonical embedding $X \hookrightarrow |\Kx|^*$. Suppose $h^{0}(X, L) = n \geq 1$, and let $s$ be a nonzero section of $L$ with divisor of zeroes $D = \sum_{i=1}^{n} k_{i}x_i$. Then in fact $\Ox(D) = L$, and the kernel of the restriction map
\[ \overline{s}^{*} \colon \Endom(L) = \Ox \to \Homom(\Ox(D), L) = \Ox \]
is the zero section of $\Ox$. Thus the scroll $\PP \uKer(\overline{s}^{*})$ disappears, and the locus $N(s)$ consists of
\[ \bigcup_{i=1}^{n} \Osc^{k_{i}-1}(\phi(X), \phi(x_{i})), \]
which is the span of the divisor $(s)$ in the usual sense.
\par
With this observation, we recover the Riemann--Kempf singularity theorem and the geometric Riemann--Roch theorem for line bundles of degree $g-1$ from Theorems \ref{gRKst} and \ref{ggRRt}.


\begin{thebibliography}{99}

\bibitem{BBCF2004} Ballico, E.; Bocci, C.; Carlini, E.; Fontanari, C.: \textsl{Osculating spaces to secant varieties}, Rend.\ Circ.\ Mat.\ Palermo, II.\ Ser.\ \textbf{53}, no.\ 3, 429--436 (2004).

\bibitem {Col2005} Coles, D.: \textsl{Vector bundles and codes on the Hermitian curve}, IEEE Trans.\ Inform.\ Theory \textbf{51} (2005), no.\ 6, 2113--2120.

\bibitem {GH1978} Griffiths, P.; Harris, J.: \textsl{Principles of Algebraic Geometry}, John Wiley and Sons Inc., USA (1978).

\bibitem {Har1977} Hartshorne, R.: \textsl{Algebraic Geometry} (Graduate Texts in Mathematics 52), Springer-Verlag Inc., New York (1977).

\bibitem {Hit2005b} Hitching, G.\ H.: \textsl{Moduli of rank 4 symplectic bundles over a curve of genus 2}, J.\ London Math.\ Soc.\ (2) \textbf{75}, 255--272 (2007).

\bibitem {HJ2008} Hitching, G.\ H.; Johnsen, T.: \textsl{Decoding of scroll codes}, ``Algebraic geometry and its applications'' (Proceedings of the first SAGA conference, Papeete, 2007, ed.\ Hirschfeld, J., Chaumine J.\ and Rolland, R.), vol.\ 5 of Number Theory and Its Applications, 294--314. World Scientific, 2008.

\bibitem {HR2004} Hwang, J.-M.; Ramanan, S.: \textsl{Hecke curves and Hitchin discriminant}, Ann.\ Sci.\ \'Ec.\ Norm.\ Sup\'er.\ (4) \textbf{37}, no.\ 5, 801--817 (2004).

\bibitem {Joh2003} Johnsen, T.: \textsl{Rank two bundles on algebraic curves and decoding of Goppa codes}, International J.\ of Pure and Applied Math., vol.\ 4, no.\ 1, 33--45 (2003).

\bibitem {Kem1983} Kempf, G.\ R.: \textsl{Abelian Integrals}, Monograf\'{\i}as del Instituto de Matem\'aticas 13, Universidad Nacional Aut\'onoma de M\'exico, Mexico (1983).

\bibitem {Kem1986} Kempf, G.\ R.: \textsl{The equations defining a curve of genus 4}, Proc.\ Amer.\ Math.\ Soc.\ \textbf{97}, no.\ 2, 219--225 (1986).

\bibitem {KS1988} Kempf, G.\ R.; Schreyer, F.-O.: \textsl{A Torelli theorem for osculating cones to the theta divisor}, Compositio Math.\ \textbf{67}, no.\ 3, 343--353 (1988).

\bibitem {LNar1983} Lange, H.; Narasimhan, M.\ S.: \textsl{Maximal subbundles of rank two vector bundles on curves}, Math.\ Ann.\ \textbf{266}, no.\ 1, 55--72 (1983).

\bibitem {LMP2006} Lanteri, A.; Mallavibarrena, R.; Piene, R.: \textsl{Inflectional loci of scrolls}, Math.\ Z.\ \textbf{258}, no.\ 3, 557--564 (2008).

\bibitem {Las1991} Laszlo, Y.: \textsl{Un th\'eor\`eme de Riemann pour les diviseurs th\^eta sur les espaces de modules de fibr\'es stables sur une courbe}, Duke Math.\ J.\ \textbf{64}, no.\ 2, 333--347 (1991).

\bibitem {Muk1995} Mukai, S.: \textsl{Vector Bundles and Brill--Noether Theory}, ``Current topics in complex algebraic geometry'', MSRI series, vol.\ 28, papers from Special Year at MSRI, Berkeley, California, 1992/1993, 145--158. Cambridge University Press, Cambridge (1995).

\bibitem {NR1969} Narasimhan, M.\ S.; Ramanan, S.: \textsl{Moduli of vector bundles on a compact Riemann surface}, Ann.\ Math.\ (2) \textbf{89}, pp. 14--51 (1969).

\bibitem {NR1975} Narasimhan, M.\ S.; Ramanan, S.: \textsl{Deformations of the moduli space of vector bundles over an algebraic curve}, Ann.\ Math.\ (2) \textbf{101}, 391--417 (1975). 

\bibitem {NR1987} Narasimhan, M.\ S.; Ramanan, S.: \textsl{$2\Theta$-linear systems on Abelian varieties}, ``Vector bundles on algebraic varieties'' (Bombay, 1984), 415--427, TIFR Stud.\ Math., 11, Tata Inst.\ Fund.\ Res., Bombay (1987).

\bibitem {NS1965} Narasimhan, M.\ S.; Seshadri, C.\ S.: \textsl{Stable and unitary vector bundles on a compact Riemann surface}, Ann.\ of Math.\ (2) \textbf{82}, 540--567 (1965). 

\bibitem {Pau2002} Pauly, C.: \textsl{Self-duality of Coble's quartic hypersurface and applications}, Michigan Math.\ J.\ \textbf{50}, no.\ 3, 551--574 (2002).

\bibitem {Pau2003} Pauly, C.: \textsl{On cubics and quartics through a canonical curve}, Ann.\ Sc.\ Norm.\ Super.\ Pisa Cl.\ Sci.\ (5) \textbf{2}, no.\ 4, 803--822 (2003).

\bibitem {Pie1976} Piene, R.: \textsl{Numerical characters of a curve in projective $n$-space}, ``Real and complex singularities'', Proc.\ Ninth Nordic Summer School/NAVF Sympos.\ Math., Oslo, 1976, 475--495. Sijthoff and Noordhoff, Alphen aan den Rijn (1977).

\bibitem {PT1990} Piene, R.; Tai, H.-s.: \textsl{A characterization of balanced rational normal scrolls in terms of their osculating spaces}, ``Enumerative geometry'', Sitges, 1987, 215--224. Lecture Notes in Math., 1436, Springer, Berlin (1990).

\end{thebibliography}
\end{document}